\documentstyle[12pt]{article}
\newcommand{\mt}[1]{\mathop{\rm #1}\nolimits}
\newcommand{\Pic}{\mathop{\rm Pic}\nolimits}
\newcommand{\Bs}{\mathop{\rm Bs}\nolimits}
\newcommand{\Sing}{\mathop{\rm Sing}\nolimits}
\newcommand{\too}{\hbox{{\rm- - -}}\to}

\font\twelvegtc=eufm10 scaled 1200

\font\ninegtc=eufm9
\font\sevengtc=eufm7

\newfam\gtcfam

\textfont\gtcfam=\twelvegtc
\scriptfont\gtcfam=\ninegtc
\scriptscriptfont\gtcfam=\sevengtc%
\font\twelveBBB=msbm10 scaled 1200
\font\tenBBB=msbm10
\font\sevenBBB=msbm7
\newfam\BBBfam
\def\BBB{\fam\BBBfam\twelveBBB}
\textfont\BBBfam=\twelveBBB
\scriptfont\BBBfam=\tenBBB
\scriptscriptfont\BBBfam=\sevenBBB

\newcommand{\CC}{{\BBB C}}

\newcommand{\ZZ}{{\BBB Z}}
\newcommand{\QQ}{{\BBB Q}}
\newcommand{\PP}{{\BBB P}}
\newcommand{\NN}{{\BBB N}}

\newcommand{\OOO}{{\cal O}}

\newcommand{\EEE}{{\cal E}}

\newcommand{\FFF}{{\cal F}}

\def\proclaim#1{\par\medskip\noindent\begingroup \it
{\bf #1.\ }}%
\def\endproclaim{\par\endgroup}

\def\remark#1{\par\medskip\noindent\begingroup \rm
{\bf #1.\ }}%
\def\endremark{\par\endgroup}

\def\definition#1{\par\medskip\noindent\begingroup \rm
{\bf #1.\ }}%
\def\enddefinition{\par\endgroup}

\def\demo#1{\par\medskip\noindent\begingroup \rm
{\it #1.\ }}%
\def\enddemo{\par\endgroup}

\def\roster{\begin{itemize}}
\def\endroster{\end{itemize}}

\title{
On the existence of good divisors \\
on Fano varieties of coindex 3}
\author{Yuri G. Prokhorov}
\date{}

\begin{document}
\maketitle
\section*{Introduction}
A normal projective variety $X$ is called {\it Fano} if some
multiplicity $-nK_X$,  $n\in\NN$ of anticanonical (Weil)
divisor $-K_X$ is an ample Cartier divisor. The number
$r(X):=\sup \{ t\in \QQ | -K_X\equiv tH, \qquad H \qquad
\mbox{{\rm is an ample Cartier divisor}} \}$
is called the {\it index} of a Fano variety $X$.
A Fano variety with  only log-terminal
singularities (see [7]) we call briefly {\it log-Fano variety},   and
a Fano variety with
only $\QQ$-factorial terminal singularities and Picard number $\rho =1$
-- $\QQ$-Fano variety.  If  $X$  is  a log-Fano variety,  then
$\Pic(X)$ is a finitely generated torsion-free group.
Therefore $r(X)\in \QQ$, $r(X)>0$.
In that case there exists the ample Cartier divisor $H$,
called {\it a fundamental divisor} of $X$,  such that
 $-K_X\equiv r(X)H$.
It is known
that $0<r(X)\le \dim(X)+1$    for any log-Fano varitiety $X$
(see e.g. [18]).
\par
We say that there exists a good divisor on $X$ if the fundamental
linear system $|H|$ is non-empty and contains a reduced irreducible  divisor
with singularities at worst the same as singularities of the variety  $X$
(e.g. if $X$ is non-singular, or has  terminal, canonical or log-terminal
singularities, then $H$ is non-singular, or has  terminal, canonical or
log-terminal  singularities, respectively).
\par
For the first time the existence of good divisors was proved in
three-dimensional non-singular case by Shokurov
 [17]. In his preprint [11]  Reid used  Kawamata's technique
for study of linear system $|H|$ and proved the existence of
good divisors for Fano threefolds with canonical singularities.
Later  this technique was applied for  Fano fourfolds of
index 2  with Picard number 1 by Wilson [14]
\footnote{For any Picard number see
[Prokhorov Yu.~G. The existence of smooth divisors on Fano fourfolds of
index 2,
{\it Russian Acad. Sci. Sb. Math.}
83 (1995) no. 1, 119--131].  (Added in translation).}
and for
log-Fano varieties $X$ of indices $r>\dim(X)-2$ by Alexeev [1].
\par
In the present paper we study the case
$r=\dim(X)-2$.  Mukai classified such non-singular Fano
varieties of any dimension under the assumption of the
existense of good divisors [9], [10] (see also [15]).
We investigate five-dimensional case and prove
the following result which is slightly weaker then Mukai's classification
claims \footnote{In the paper [Prokhorov Yu.~G.
On the existence of good divisors on Fano varieties  of  coindex
3, II
{\it  Contemporary Math. and its Appl. Plenum. 24 (1995)} (to appear)]
the author generalized this result and proved that $H$ is non-singular
for $\dim X=5$ and non-singular $X$.
(Added in translation).}.

\proclaim{Theorem 1} Let $X$ be a non-singular Fano fivefold
of index 3 and let $H$ be a fundamental divisor on $X$.  Then the
linear system  $|H|$ contains an irreducible divisor with only canonical
Gorenstein singularities.
\endproclaim
The main idea of the proof of the Theorem is  to investigate the singular
locus $\Sing(H)$  of a general divisor $H$.  Here the "bad" case is
the case when $\dim \Sing(H)=3$. In this situation we study
a three-dimensional component of $\Sing(H)$, applying
the minimal model program. In fact using similar arguments one can prove
more general result:
\proclaim{Theorem 2} Let $X$ be a $n$-dimensional log-Fano variety
of index $n-2$ and let $H$ be a fundamental divisor on $X$.
Assume that in dimension $n-2$ flip-conjectures I and II hold  (see [7]).
Then if the linear system $|H|$ is not empty and has no fixed
components, then it contains an irreducible divisor with only log-terminal
singularities.
\endproclaim
\remark{Remark} If  in notations of Theorem 0.2  $X$  has only canonical
Gorenstein singularities, then $|H|\ne \emptyset$  (see 1.1 ).
\endremark
\subsection*{Notations and conventions:}
\par
\roster
\item[$\Bs |H|$:] the scheme-theoretic base locus of a linear system $|H|$;
\item[$\Sing (X)$:] the singular locus of a variety $X$;
\item[$\equiv$:] numerical equivalence of cycles;
\endroster

The ground field is assumed to be the field of complex numbers $\CC$.
 We will use the basic definitions and concepts of the minimal
model program (see [7]).

\section*{1. Preliminary results}
{}From Kawamata-Viehweg vanishing theorem for varieties with log-terminal
singularities (see [7]) we get
\proclaim{1.1. Lemma   [1]} Let $X$ be a  $n$-dimensional
log-Fano variety  of index $r$,
$H$ be a fundmental divisor on $X$ and $H^n=d$. Then
\roster
\item[(i)]  If  $r>n-2$, then $\dim|H|=n-2+d(r-n+3)/2>0$. Moreover
$r=n-3+2k/d$ for some $k\in\NN$, $k>d/2$;
\item[(ii)] If $r=n-2$ and $X$  has only canonical Gorenstein singularities,
then $\dim|H|=g+n-2>0$, where $2g-2=d$, $g\in \ZZ$, $g\ge 2$.
\endroster
\endproclaim

In  [2] Fujita defined $\Delta$-genus of a polarized $n$-dimensional
variety $(X,H)$ as
 $\Delta (X,H)=H^n -h^0 (X,\OOO_X (H))+n$.
\proclaim{1.2. Corollary}  In notations (i) of Lemma 1.3 we have
$\Delta(X,H)=d-k+1$.
\endproclaim
\proclaim{1.3. Theorem  [2]}  $\Delta (X,H)\ge\dim(Bs|H|)+1\ge 0$   for any
polarized variety $(X,H)$. Moreover if
$\Delta (X,H)=0$, then the divisor $H$ is very ample.
\endproclaim
The following is a consequense from the classification of polarized
varieties of $\Delta$-genus zero [2]
\proclaim{1.4. Corollary  (Fujita)} Let $X$ be a $\QQ$-Fano threefold
of index $r>2$. Then either
\roster
\item[] $X=\PP^3$,
\item[] $X=Q\subset \PP^4$ is a smooth quadric, or
\item[] $X=X_4\subset \PP^6$  is
a projective cone over the Veronese surface.
\endroster
\endproclaim
\proclaim{1.5. Theorem  [4], [12]} Let $X$ be a three-dimensional Fano
variety of index 2
with only  canonical  Gorenstein singularities, let
$$ be a fundamental divisor on $X$ and $d:=H^3$. Then
\roster
\item[(i)] if  $d\ge 3$, then $H$ is very ample,
\item[(ii)] if $d=2$,  then the linear system  $|H|$ defines a finite
morphism $X\to \PP^3$ of degree 2.
\endroster
\endproclaim

\section*{2. Kawamata's technique}
We describe  brefly the technique of resulution of base loci
of linear systems on Fano varieties in connection with our
situation ( see [11], and also  [1],[14]).
\par
Let $X$ be a non-singular Fano fivefold of index 3 and let $H$ be a
fundamental divisor on $X$. By Lemma 1.1, $|H|\ne \emptyset$.
Then there exist a resolution $f:Y\to X$ and a divisor $\sum E_i$ on
 $Y$ with only simple normal crossings such that
$$
1)\qquad K_Y\equiv f^* K_X +\sum a_i E_i,\quad a_i\in \ZZ,\quad a_i\ge 0,
\eqno(2.1)
$$
where $a_i\ne 0$ only if $f_*E_i =0$;

$$
2)\qquad f^* |H|=|L|+\sum r_i E_i,\quad  r_i \in \ZZ,\quad  r_i\ge 0,\eqno(2.2)
$$
where the linear system  $|L|$ is free;
\par\noindent
3) $\QQ$-divisor $qf^* H-\sum p_i E_i$  is ample for some
$0<p_i\ll 1$, $0<q<\min \{ 1/r_i | r_i\ne 0\}$.
\par
Set $c=\min\{ (a_i +1-p_i )/r_i  | r_i\ne 0 \}$. By changing coefficients
$p_i$ a little one  can attain that the minimum will achive for only
one index, say for $i=0$. By Kleiman's criterion for ampleness, the
following  $\QQ$-divisor
$$
N=N(t)=tf^*H+\sum(-cr_i+a_i-p_i)E_i-K_Y\equiv
$$
$$
\equiv cL+f^*(t-c+3)H-\sum p_i E_i,\quad t\in\ZZ
$$
is ample for $t-c+3\ge q>0$. Since $-cr_0+a_0-p_0=-1$  and
$-cr_i +a_i -p_i >-1$ for $i\ne 0$,  for the rounding-up  of $N$
we have $N =tf H-K +A-E$,  where  $E=E_0$,  $A\ge 0$ and  $A$  consists of
exceptional for $f$ components of $\sum E_i$.
\proclaim{2.1. Lemma   [1],[11]} If $4\ge c+q$, then
$H^0(E,\OOO_E(f^*H+A))=0$.
In particular, $H^0(E,\OOO_E(f^*H))=0$.
\endproclaim
\demo{Proof} By the Kawamata-Viehweg vanishing theorem (see [6],[13])
the following sequence
$$
0\to H^0(Y,\OOO_Y(f^*H+A-E))\to H^0(Y,\OOO_Y(f^*H+A))\to
H^0(E,\OOO_E(f^*H+A))\to 0
$$
is exact.
If $H^0(E,\OOO_E(f^*H+A))\ne 0$, then we get a contradiction
with the fact that  $E$ is a fixed component of the linear system
 $|f^*H+A|$.
\enddemo
\proclaim{2.2. Lemma [1],[11]} For constants $a_i$  and $r_i$ in formulas
 (2.1) and (2.2) the inequality $a_i +1\ge r_i$ holds for any $i$.
\endproclaim
\proclaim{2.3. Lemma} Assume that $a_i +1=r_i$ for some  $i$. Then
then   there are the following possibilities
for the divisor $E=E_0$:
\roster
\item[(i)] $a_0 =0$, $r_0 =1$, $\dim(f(E))=4$, and $f(E)$ is a fixed component
of  $|H|$ of multiplicity $1$;
\item[(ii)] $a_0=1$, $r_0=2$, $\dim(f(E))=3$, and a general divisor $|H|$ has
only
double singularities along $f(E)$.
\endroster
\endproclaim
\demo{Proof}
Since $a_i+1=r_i$, we have $c\le1$.  Set  $d:=\dim(f(E))$.
Consider the following polinomial of degree $\le d$:
$p(t)=\chi (\OOO_E(tf^*H+A))$. For $t>q-2$  $\QQ$-divisor $N(t)$
is ample and by the Kawamata-Viehweg vanishing theorem (see [6],[13]),
$p(t)=h^0(E,\OOO_E(tf^*H+A))$,  i.~e.
the polinomial $p(t)$ has zero for $t=-1$, and also, by Lemma 2.1,
one more zero for $t=1$. On the other hand, $p(0)=h^0(E,\OOO_E(A))=1$.
Hence,  if  $d\le 2$, then $p(t)\le 0$  for   $t\gg 0$,  a contradiction
with ampleness of $H$. Therefore, $d=4$ and $d=3$.
\enddemo
\proclaim{2.4. Lemma} Let $E_k$ be a component of divisor $\sum E_i$
such that $\dim f(E_k)=3$ and $x\in f(E_k)$ be a general point. Then
a general divisor $H$ is non-singular or has only double normal
crossing singularities\footnote{
i.~e. $x\in H$ is analytically isomorphic to $(0\in\{ x_1x_2=0\} )\subset
(0\in\CC^{5})$. (Added in translation).} at $x$.
\endproclaim
\demo{Proof} A general surface $X'=D_1\cap D_2\cap D_3$, where $D_i\in |mH|$,
$m\gg 0$ is non-singular.
Moreover for the corresponding resolution
$$
\begin{array}{ccc}
Y'&\hookrightarrow&Y\\
\downarrow\lefteqn{{\scriptstyle f'}}&&\downarrow\lefteqn{{\scriptstyle f}}\\
X'&\hookrightarrow&X\\
\end{array}
$$
we have
$$
K_Y=f^*K_{X'}+\sum a_iE_{ij}',\qquad\qquad f^*|H'|=|L'|+\sum r_iE_{ij}'
$$
where $E'_{ij}$ is a component of $Y'\cap E_i$. Assume that a general divisor
 $H$ has  worse singurarity at $x$ than double normal crossings. Then
the curve $H'$ has  worse singurarity at the  point $x$ than ordinary  double.
If  $x\in H'$ is not simple cusp,  then for the "first" blow-up at
$x$ we have $a_1=1$, $r_1\ge 2$, and for the "second" blow-up ---  $a_2=2$,
$r_2\ge 4$, a contradiction with Lemma  2.2. But if  $x\in H'$ is a simple
cusp, then
similar we have
$a_1=1$,  $r_1=2$,  $a_2 =2$,  $r_2=3$,   $a_3=4$,   $r_4=6$,  again
it contradicts to
Lemma  2.2.  Thus  $H'$ has at  $x$ only ordinary double
point.
This proves our lemma.
\enddemo

\section*{3. Some corollaries from the minimal model program}
\proclaim{3.1. Theorem    [8]}
Let $S$ be a projective three-dimensional variety.
Then  there exist birational morphisms $\mu:W'\to S$,
$\upsilon:W\to W'$ and $\tau=\mu\circ\upsilon:W\to S$ with the following
properties:
\roster
\item[(i)] $W$ has only terminal $\QQ$-factorial singularities,  $W'$
has only canonical singularities;
\item[(ii)]  $K_W$  is   $\tau$-numerically effective, $K_{W'}$ is $\mu$-ample;
\item[(iii)] morphism  $\upsilon$ is crepant, i.~e. $K_W=\upsilon^*K_{W'}$;
\item[(iv)] if $\alpha:V\to S$ is any birational morphism,  $V$ has only
terminal
$\QQ$-factorial singularities, then the map
$\tau^{-1}\circ\alpha: V\too W$ is a composition of contractions of
extremal rays and flips.
\endroster
\endproclaim
\definition{Definition}
Varieties  $W$ and $W'$ from Theorem  3.1 are called
 by {\it terminal} and {\it canonical modifications} of $S$, respectively.
\enddefinition
\proclaim{3.2. Theorem    [3,  \S 3-4]} Let  $W$ be a three-dimensional
projective variety with only $\QQ$-factorial terminal singularities
and $N$ be a numerically effective big Cartier divisor on $W$.
Assume that there exists an extremal ray on $W$ is generated by a class
of a curve $\ell$ such that $(K_W+2N)\cdot\ell<0$ and let
$\varphi:W\to C$ be its contraction. Then one of the following holds:
\roster
\item[(i)] $\varphi :W\to C$ is a birational morphism and $\ell\cdot N=0$;
\item[(ii)]  $W$ and $C$ are smooth, $\dim C=1$, $W=\PP_C(\EEE)$ and
$N=\OOO_{\PP_C(\EEE)}(1)$,  where $\EEE$ is a locally free sheaf
of rank 3 on $C$;
\item[(iii)]  $C$ is a point  and  $W$ is a $\QQ$-Fano variety of
index $r>2$.
\endroster\endproclaim

\proclaim{3.3. Proposition} Let $\alpha:\hat{S}\to S$ be a birational
morphism of three-dimensional normal projective varieties, where   $\hat{S}$
has only $\QQ$-factorial terminal singularities, and let  $M$ be an ample
Cartier divisor on $S$. Assume that $K_{\hat{S}}=-2\hat{M}-\hat{R}+\hat{D}$,
where $\hat{M}=\alpha^*M$, $\hat{R}$ is an effective $\QQ$-divisor,
$\hat{D}$ is a $\QQ$-divisor such that all the components of $\hat{D}$
 are contracted by the morphism $\alpha$.  Then for terminal and canonical
modifications $\tau:W\to S$ and  $\mu:W'\to S$, and also for Cartier divisors
$N=\tau^*M$ and  $N'=\mu^*M$ we have only one of the following possibilities:
\roster
\item[(i)]  $W=\PP_C(\EEE)$, $N=\OOO_{\PP(\EEE)}(1)$, where $C$ is a
smooth curve, $\EEE$ is a locally free sheaf on   $C$ rank 3;
\item[(ii)]
 $S=W=W'=\PP^3$, $M=N=N'=\OOO_{\PP^3}(1)$;
\item[(iii)]
 $S=W=W'=Q\subset \PP^4$ is a non-singular quadric,
$M=N=N'=\OOO_Q(1)$;
\item[(iv)]
  $S=W=W'=W_4\subset\PP^6$ is a cone over the Veronese surface,
$M=N=N'=\OOO_W(1)$;
\item[(v)] $S=W'$ is a Fano variety of index $2$ with only canonical Gornstein
singularities,  $-K_{W'}=2N'=2M$,  $W$  has only isolated  cDV-points,
 $-K_W=2N$ is a numerically effective big divisor.
\item[(vi)] $S=W=W'=\PP^3$ , $M=N=N'=\OOO_{\PP^3}(2)$.
\endroster\endproclaim
\demo{Proof}  By Theorem  3.1, the birational map
$\beta:=\tau^{-1}\circ\alpha:\hat{S}\too W$  is a composition of contractions
of extremal rays and flips. In particular the inverse map
$W\too\hat{S}$ doesn't contract divisors. Therefore
$$
K_W=-2N-B+G,                   \eqno (3.1)
$$
where $N:=\tau^*M$, $B$ is an effective $\QQ$-divisor on $W$, and
 $G$ is a $\QQ$-divisor such that all its components
 are contracted by the morphism $\tau$. From  (3.1) we have
$$
K_{W'}=-2N'-B'+G',\eqno (3.2)
$$
where $N':=\tau^*M$, $B'=\upsilon_*B$ is an effective $\QQ$-Weil divisor
on $W$ and $G'=\upsilon_*G$ is $\QQ$-Weil divisor such that all the components
of $G$
 are contracted by the morphism $\mu$.
First we assume that $B=G=0$. Then  $B'=G'=0$ and $-K_{W'}=2N'$.
If $c$ is a curve in a fiber of $\mu$, then
$K_{W'}\cdot c=-2N'\cdot c=0$, this  contradicts to  $\mu$-ampleness of
$K_{W'}$. So $\mu$ is a finite birational morphism  on the  normal variety
$S$. Hence  $\mu$ is an isomorphism and   $-K_{W'}=2N'$ is an ample divisor.
We obtain cases (v), (vi). Now let  $B\ne 0$ or $G\ne 0$.
We claim that $K_W+2N$ is not numerically effective.
Indeed, in the opposite case $-B+G=K_W+2N$ is numerically effective
and for any irreducible curve $á$ on $S$ we have $0\le (-B+G)\cdot \tau^*c
=-B\cdot \tau^*c$. Thus $\tau_*B=0$, i.~e. we may assume that $B=0$.
But then by the Base Point Free Theorem  (see e.~g. [18],[7]),
 the lineear system
$|kG|$ is free for $k\gg 0$, it contradicts to contractedness of
the divisor $G$. Futher by the Cone Theorem ( see [7]), there exists
an extremal ray on $ W$ generated by a class of a curve $\ell$
such that   $(K_W+2N)\cdot\ell<0$. Let $\varphi:W\to C$ be the contraction of
this extremal ray.
Note that if $\ell\cdot N=0$, then $\ell$ is contained in a fiber of morphism
$\tau$. It contradicts $\tau$-numerical effectiveness of the divisor $K_W$.
Using Theorem   3.2  and Corollary 1.7, we obtain cases  (i)-(iv).
\enddemo
\proclaim{3.4. Lemma} In the case  (i) of Proposition 3.3  the curve $C$
can be rational or elliptic.
\endproclaim
\demo{Proof} Denote by  $F$ the class of a fiber of the projection
 $\varphi:\PP(\EEE)\to C$. We have the standard formula
$$
K_W=-3N+\varphi(\det(\EEE)+K_C).
$$
Set $d:=\deg\EEE$. Then
$$
K_W\equiv -3N+(2g(C)-2+d)F.     \eqno              (3.3)
$$
Comparing (3.1) and  (3.3) we obtain
$$
B\equiv N-(2g(C)-2+d)F+G     \eqno (3.4)
$$
We may assume that the morphism $\tau$ does not contract
components of $B$. In the Chow ring $A(W)$ the following conditions are
satisfied
$$
N^3=d,\qquad N^2\cdot F=1,\qquad F^2=0      \eqno (3.5)
$$
Moreover   $d=N^3=M^3\ge 1$ because  $N=\tau^*M$.
\par
Assume that an irreducible  divisor $P$ is contracted by
the morphism  $\tau$.  Then $P\equiv aN+bF$,  $a,b\in \ZZ$,  and from
 (3.5)  one can see $0=N^2\cdot P=da+b$ and  $0\le N\cdot F\cdot P=da$,
i.~e.  $P\equiv a(N-dF)$, $a\in \NN$. If $\tau$ contracts one more
irreducible divisor $P'$, then $ P'\equiv a'(N-dF)$,  $a'\in \NN$  and
$0\le P\cdot P'\cdot N=-aa'd<0$. The contradiction shows that $\tau$
may contracts at most one irreducible divisor on $W$, i.~e.  $G=0$ or  $G=pP$,
where $p\in \QQ$ and $P\equiv a(N-dF)$ is an irreducible exceptional
  divisor.
First assume that $B=0$. Then from  (3.4) and  (3.5) we have
$N\equiv (2g(C)-2+d)F-G$, $d=N^3=(2g(C)-2+d)N^2\cdot F=2g(C)-2+d$,
i.~e.   $g(C)=1$.
But if $B\ne 0$, then $0<N^2\cdot B=2-2g(C)$, i.~e.  $g(C)=0$.
\enddemo
\proclaim{3.5. Corollary}  In conditions of  Proposition 3.3  we have
$H^0(S,\OOO_S(M))\ne 0$.
\endproclaim
\demo{Proof} Since $S$ is normal,  then  it is sufficient to prove only that
$H^0(W,\OOO_W(N))\ne 0$  or  $H (W',\OOO_{W'}(N'))\ne 0$.
The last two inequalities are easy in cases (ii), (iii), (iv) and (vi),
and in the case (v) they follow from Lemma 1.1. Consider the case (i). Then
$H^0(W,\OOO(N))=H^0(C,\EEE)$ and  by Riemann-Roch on $C$ we have
$h^0(C,\EEE)=3(1-g)+d+h^1(C,\EEE)\ge 1$.
\enddemo

\section*{4. Proof of Theorem 1}
     Let $X$ be a non-singular Fano fivefold of index 3 and
let $H$ be a fundamental divisor on $X$. Then by Lemma 1.1, $\dim|H|\ge 6$.
It follows from results of [16] that $\Pic(X)\simeq\ZZ\cdot H$, except
the following three cases: $X=\PP^2\times \QQ^3$,
$X$ is a divisor of bidegree (1,1) on $\PP^3\times \PP^3$ or
$X$ is the blow-up  of $\PP^5$ along $\PP^1$. In all of these cases
it is easy to check directly that the linear system  $|H|$ contains
a smooth divisor. Thus we will suppose that  $\Pic(X)\simeq\ZZ\cdot H$.
In particular we suppose that the linear system $|H|$ has no fixed components.
We will use
all the notations of Section 2. The following proposition
is the main step in the proof of Theorem 1.
\proclaim{4.1. Proposition} $r_i <a_i +1, \forall i$.
\endproclaim
\demo{Proof} According to Lemma 2.2,  $r_i\le a_i +1$. Assume  that
$r_i =a_i +1$  for some $i$.  Then $ c:=\min  \{ (a_i +1-p_i )/r_i   |
r_i\ne 0\} \le 1$. Set  $Z:=f(E)$ (with reduced structure).  By Lemma 2.3
a general divisor $H$ has only double normal crossings singularities in a
general
point  $x\in Z$. We may suppose that the resolution $f:Y\to X$
is a composition  $h:Y\to\tilde{X}$ and $g:\tilde{X}\to X$, where
$g:\tilde{X}\to X$
is a resolution of singularities of a general divisor $H$
(i~e. we fix some general divisor $H$ and suppose that $g$ is a composition
of blow-ups with centers in subvarieties contained in singular locus of
$H$). We have
$$
K_{\tilde{X}}=g^*K_X+\sum a_i'E_i',      \eqno       (4.1)
$$
$$
g^*H=\tilde{H}+\sum r_i' E_i',\qquad    r_i\ge 2,     \eqno   (4.2)
$$
(here $\tilde{H}$ is a proper transform of $H$, $\sum E'_i$ is   an
exceptional divisor).
For corresponding  $E_i$, $E'_i$;\quad $r_i$, $r'_i$  and $a_i$, $a'_i$
the following conditions are satisfied:
$$
h(E_i)=E_i',\qquad  r_i\le r'_i,\qquad  a_i =a'_i.
$$
Denote by $E'$ the component of the divisor $\sum E'_i$, corresponding $E$.
According to Lemma 2.4, there are two possibilities:
$$
(I)\qquad\qquad
\tilde{H}\cap E'=\tilde{Z}+\tilde{Z_1}+\tilde{Z_0},   \eqno         (4.3)
$$
where restrictions $g:\tilde{Z}\to Z$ and $g:\tilde{Z_1}\to Z$ are
birational, $g_*\tilde{Z_0}=0$;
$$
(II)\qquad\qquad
\tilde{H}\cap E'=\tilde{Z}+\tilde{Z_0},   \eqno        (4.4)
$$
where the restriction $g:\tilde{Z}\to Z$ is generically finite of degree 2,
$g_*\tilde{Z_0}=0$.
\par
We study the variety   $\tilde{Z}$ and corresponding morphism
$\psi:=g|_{\tilde{Z}}:\tilde{Z}\to Z$. Denote by $T$ the
Cartier  divisor $H|_{Z}$  on $Z$.
\enddemo
\proclaim{4.2. Lemma}
$K_{\tilde{Z}}=-2\tilde{T}-\tilde{B}+\tilde{G}$,
where $\tilde{T}$, $\tilde{B}$, $\tilde{G}$ Cartier divisors on $\tilde{Z}$,
$\tilde{G}$ is contracted by the morphism $\psi$, $\tilde{B}$
is effective and  $\tilde{T}:=\psi^*T$.
\endproclaim
\demo{Proof} Consider for example the case I. The variety $\tilde{Z}$
is Gorenstein because it is divisor on a non-singular variety. By the
ajunction formula we have
$$
K_{\tilde{Z}}=(K_{\tilde{H}} +\tilde{Z})|_{\tilde{Z}}=
(K_{\tilde{X}}|_{\tilde{H}}+\tilde{H}|_{\tilde{H}}+\tilde{Z})|_{\tilde{Z}}
$$
$$
=-2\tilde{T}+\left.\left(\sum (a'_i-r'_i)E'_i|_{\tilde{H}} +E'|_{\tilde{H}}-
\tilde{Z_1} -\tilde{Z_0}\right)\right|_{\tilde{Z}}=
$$
$$
-2\tilde{T}+\left.\left(\sum_{i\ne 0}(a'_i-r'_i)E'_i|_{\tilde{H}}
-\tilde{Z_1} -\tilde{Z_0}\right)\right|_{\tilde{Z}}
$$
(because  $a'_0=1$, $r'_0=2$). It follows from 2.6 that $\tilde{Z}$
is not contained in $E'_i$ for $i\ne 0$.
It is sufficient to prove  that in the last formula $a_i'\le r_i'$,
if the corresponding divisor $E'_i|_{\tilde{Z}}$ is not
contracted by the morphism $\psi$. If
$\dim(g(E'_i))\le 1$, then for every component $F'$ of the divisor
$E'_i|_{\tilde{Z}}$ we have
$\dim(\psi(F'))\le \dim(g(E'))\le 1<\dim(F')$,  i.~e.  $F'$ is contracted
by the morphism $\psi$ is this case.
But if $\dim(g(E'_i))=3$,  then  by  Lemma   2.4,  $r_i'\ge a_i'$.
It remains to prove only that the case $\dim(g(E'_i))=2$ and $a_i'>r_i'$
is impossible.
Suppose that $a_i'>r_i'$ for some $i$   and consider the following
partial ordering of the set
$\{ E_j' \}$:\qquad  $E'_j\succ E'_k$ (or simply $j\succ k$), if
and only if
the center of the $k$-th blow-up is contained in  the $j$-th
exceptional  divisor.
Then the divisor $E_i'$ cannot be maximal, otherwise  $a_i'=2$,
$r_i'\ge 2$, that contradicts  our assumption.
Let $\sigma_i:X_i\to X_{i-1}$ be the blow-up of a smooth surface
corresponding to  the divisor  $E_i'$.
Then  $K_{X_i}=\sigma_i^*K_{X_{i-1}}+2\bar{E_i}$,
$\sigma_i^*H_{i-1}=H_{i} +m\bar{E_i}$,   $m\ge 2$,  where   $\bar{E_i}$
is the exceptional divisor of the blow-up  $\sigma_i$  (its proper transform on
$\tilde{X}$ is the divisor $E '_i$), $H_{i-1}$    and $H_{i}$ are proper
transforms of $H$ on
$X_{i-1}$
    and
$X_i$ , respectively. We obtain
$$
a'_i=2+\sum_{j\succ i} a'_j ,\qquad
r'_i=m+\sum_{j\succ i} r'_j ,\qquad m\ge 2.
$$
Hence   $\sum_{j\succ i} a'_j  >\sum_{j\succ i} r'_j$
(we assume,   that   $a'_i>r'_i$).
Therefore there exists $j$ such  that $j\succ i$,  $\dim(g(E'_j))\ge 2$ and
$a'_j>r'_j$.  Thus we may obtain the maximal element
 $E_j'$.  According  the above for  that maximal element we have
$\dim(g(E'_j))\ne 2$ , i.~e. $\dim(g(E'_j))=3$, but then by  Lemma  2.4
$a'_j=1$,
$r'_j\ge 2$ and the inequality $a'_j>r'_j$ is impossible. Lemma is proved.
\enddemo
\proclaim{4.3. Lemma}
There exists a resolution
$\sigma:\hat{Z}\to  \tilde{Z}$  such that
$$
K_{\hat{Z}} =-2\hat{T}-\hat{B}+\hat{G},\eqno     (4.6)
$$
where  $\hat{T}$, $\hat{B}$, $\hat{G}$ are $\QQ$-divisors on $\hat{Z}$,
$\hat{G}$ is contracted by the morphism $\hat{\psi}:=\psi\circ\sigma$,
$\hat{B}$ is effective and $\hat{T}=\hat{\psi}^*T$.
\endproclaim
\demo{Proof}
It is sufficient to prove that formulas  (4.5), (4.6)
are preserved under one blow-up $\sigma_1:\hat{H}\to\tilde{H}$
of non-singular $k$-dimensional subvariety in $\tilde{Z}$ ($k\le 2$). Let
$\hat{E}$ be an exceptional divisor of $\sigma_1$ and
$\hat{Z}$ be the proper transform of $\tilde{Z}$. Then
$K_{\hat{H}}=\sigma_1^*K_{\tilde{H}}+(3-k)\hat{E}$,
$\hat{Z}\sim\sigma_1^*\tilde{Z}-p\hat{E}$, $p\ge 2$, $p\in \ZZ$
(remind that we assumed that $\tilde{Z}$  is singular along our suvariety).
Thus by the ajunction formula
$K_{\hat{Z}}=\sigma_1^*K_{\tilde{Z}}+(3-k-p)\hat{E}|_{\hat{Z}}$,
where  $\hat{E}$ is contracted by the morphism
$\sigma_1$  (the case $k\le 1$) or $k=2$, $3-k-p=1-p<0$, i.~e.
formulas  (4.5), (4.6) are preserved .
\enddemo
Now applying Proposition 3.3 we obtain the diagram below.
%%DIAGRAM
$$
\begin{array}{cccccc}
\hat{Z}=&\hat{S}&\stackrel{\beta}{{\too}}&W'&\stackrel{\upsilon}{\longrightarrow}&W\\
        &
&\stackrel{\alpha}{\searrow}&\downarrow&\stackrel{\tau}{\swarrow}&\\
        &       &                        &S &                            & \\
        &       &                        &\downarrow\lefteqn{\scriptstyle{\pi}}
&                            & \\
        &       &                        &Z' &                            & \\
        &       &                        &\downarrow\lefteqn{\scriptstyle{\nu}}
&                            & \\
        &       &                        & Z &                            & \\
\end{array}
\eqno (4.8)
$$
where $\nu :Z'\to Z$ is the normalization of $Z$, $\pi:S\to Z'$ is
the normalization of $Z'$ in the function field of $Z$.
In the case I $\pi$ is an isomorphism. In the case II
$\pi$ is a finite morphism of degree 2, and $Z'\simeq S/\Gamma$,
where $\Gamma$ is the group of order 2.
The other notations we fix the same as in
Theorem 3.1. For $W$  and $W'$ we have only posiibilities
  (i)-(vi)  from Proposition 3.4.
Set $M:=\pi^*T'$.
\par
The divisor $E$ is a proper transform of the exceptional divisor
 of some blow-up of a non-singular model $Z''$ of the  variety $Z$ on $X''$.
Let $Z''\to Z$ be a corresponding birational morphism.
Clearly that it factors through the normalization
$Z''\stackrel{\nu'}{\longrightarrow} Z'\stackrel{\nu}{\longrightarrow} Z$.
If $r_i\ge a_i +1$ for some $i$, then by Lemma  2.1
$H^0 (E,\OOO_E (f^*H))=0$.
Therefore to prove Proposition 4.1 it is sufficient to show
 that
$$
H^0(Z'',\nu^{\prime *} \OOO_{Z'}(T'))\ne 0
$$
or
$$
H^0(Z',\OOO_{Z'}(T'))\ne 0,   \eqno    (4.8)
$$
where $T':=\nu^*T$,\quad $T:=H|_Z$. In the case I the inequality  (4.8) follows
from
3.5. Consider the case II. Then  $Z'=S/G$ and
$H^0(Z',\OOO_{Z'}(T'))=H^0(Z',\pi_*\OOO_S\otimes\OOO_{Z'}(T'))^\Gamma=
H^0(S,\OOO_S(M))^\Gamma$.
Further we assume that $H^0(S,\OOO_S(M))^\Gamma=0$.

\proclaim{4.4. Lemma} If  $H^0(S,\OOO_S(M))^\Gamma=0$, then the rational map
$S\too S_0\subset \PP^{\dim |M|}$ assosiated with the linear system
$|M|$  is not birational on its image.
\endproclaim
\demo{Proof}  If $H^0(S,\OOO_S(M))^{\Gamma}=0$, then the action of
 $\Gamma$  on the linear system  $|M|$ is trivial and if $|M|$
defines the rational map $S\too S_0$, then
it factors  through quotient morphism by  $\Gamma$:\quad
$S\stackrel{\pi}{\longrightarrow}S/\Gamma=Z'\too S_0$ .
\enddemo
    Lemma  4.4 excludes cases (ii),(iii),(iv),(vi) of Proposition
3.3. In the case (v) of Proposition 3.3 by Theorem 1.5, possibility  $M^3\ge 3$
is also
impossible. On the other hand, we have  $M^3=(\pi^*T')^3=2(T')^3\ge 2$.
If $M^3=2$, then again by Theorem 1.5  $S_0=\PP^3$ and  $S\too S_0$
is a finite morphism of degree  2, so  $S/\Gamma=S_0=\PP^3$.
It gives us inequality (4.8).
\par
Finally let $S$ and $W$ be as in  (i) of Proposition 3.3.
By Lemma 3.4, $g(C)=1$ or   g(C)=0.  If g(C)=0,  then  the sheaf
$\EEE$ on $C=\PP^1$ is decomposiable:
$\EEE=\OOO(d_1)\oplus\OOO(d_2)\oplus\OOO(d_3)$, where $d_i\ge 0$.
Then the linear system $|N|=|\OOO_{\PP(\EEE)}(1)|$ defines
 a birational map (see e.~g.  [5]), again we have a contradiction with
Lemma 4.4.  The last case 3.3  (i)  and
$g(C)=1$ is treated  similar to  3.3  (v).
Here insted Theorem 1.5 we use the following
\proclaim{4.5. Lemma}  Let  $W$  and $N$ be as in (i) of Proposition 3.3.
Assume  also that $g(C)=1$ and $d:=N^3\ge 2$.
Then the linear system $|N|$ on   $W$  defines  a
morphism  $W\to W_0\subset \PP^{d+1}$.
Moreover this morphism is finite of degree 2 if $d=2$ and
birational if $d\ge 3$.
\endproclaim
\demo{Proof} It follows from proof of Lemma 3.4 that
the morphism $\tau$ contracts a surface $P\equiv a(N-dF)$,
$a\in \NN$ to a curve. We claim that $a=1$.  Indeed in the opposite case
every divisor from  $|N|$ is irreducible. Then $N$ is a smooth
geometrically ruled surface over $C$  and  $P|_{N}$  is an effective
 divisor  with negative self-intersection number. It is known in this case
that $P|_{N}$    is a section of the ruled surface $N$,
hence $a=1$. Therefore $P$ is also a non-singular geometrically
ruled surface over $C$. The morphism $\tau$ maps $P$ onto curve
and fibers of $P\to C$ are not contracted. It gives us  that
$P\simeq C\times \PP^1$. Consider the exact sequence
$$
0\to H^0 (\varphi^*\FFF)\to H^0(\OOO_W(N))\to H^0 (\OOO_P(N))\to 0,\eqno
(4.9)
$$
where $\FFF$ is a sheaf of degree $d$ on $C$ such that
$\OOO_W(N-P)=\varphi^*\FFF$.  It is easy to see that the sheaf
$\OOO_P(N)$ on  $P\simeq C\times \PP^1$   has bidegree (0,1). So
$h^0(\OOO_P(N))=2$,  $h^0(\OOO_W(N))=d+2$ and   $P$
is not a fixed component in $|N|$. From (4.9) and $d\ge 2$ we get, that
the linear system $|N|$ is free and defines a morphism
$\tau_0:W\to W_0\subset\PP^{d+1}$. Moreover  $\deg\tau_0\cdot \deg  W_0
=N^3=d$.
Applying  to variety $W_0 \subset \PP^{d+1}$ the inequality
$\deg W_0 \ge \mt{codim} W_0 +1$  we obtain $\deg \tau_0 =2$ for $d=2$, and
$\deg \tau_0  =1$  for $d\ge 3$.
\par
Thus Proposition  4.1 is proved.
\enddemo
Now we prove that Proposition 4.1 implies Theorem 1.
By our assumption every divisor $H$ is irreducible and for
 general $H$ we have $\dim\Sing(H)<3$ (see 2.4, 4.1).
Therefore a general divisor  $H$ is normal.
For such  $H$ by the ajunction
formula and by Proposition 4.1 we have
$K_L =f|_L^* (K_H)+\sum (a_i-r_i)E_i|_L$  and  $a_i-r_i\ge 0$,
where  $f|_L:L\to H$  is the corresponding resolution of singularities of
$H$. It is equivalent the fact that  $H$ has only canonical singularities.
This proves our theorem.

\par\bigskip\noindent
Chair of Algebra \\
Department of Mathematics and Mechanics \\
Moscow State University \\
Moscow, 117 234, Russia
\par\medskip\noindent
\begin{tabular}{ll}
E-mail:\quad& prokhoro@mech.math.msu.su
 \end{tabular}

\end{document}